\documentclass[12pt]{article}
\usepackage{amsmath,amsthm,amssymb}
\usepackage{mathtools}

\newtheorem{theorem}{Theorem}[section] 
\newtheorem{lemma}{Lemma}[section]
\newtheorem{cor}{Corollary}[section]
\newtheorem{remark}{Remark}[section]

\newtheorem{example}{Example}[section]
\newtheorem{ass}{Assumption}[section]

\newtheorem{proposition}{Proposition}[section]

\numberwithin{equation}{section}

\newcommand{\D}{{\rm d}}
\newcommand{\dx}{\, \D x}

\newcommand{\du}{\, \D u}

\newcommand{\rz}{\mathbb{R}}

\newcommand{\iom}{\int_{\Omega}}

\newcommand{\klauf}{\left(\begin{array}}
\newcommand{\klzu}{\end{array}\right)}

\title{Some geometric properties of nonparametric $\mu$-surfaces in $\rz^3$}
\author{Michael Bildhauer \& Martin Fuchs}
\date{}

\newcommand{\reff}[1]{(\ref{#1})}

\usepackage{hyperref}

\hypersetup{
   pdfnewwindow=true,
   colorlinks=false,
   linkcolor=black,
   citecolor=green,
   filecolor=magenta,
   urlcolor=cyan,
}

\begin{document}

\parindent0em
\maketitle

\newcommand{\op}[1]{\operatorname{#1}}
\newcommand{\bv}{\op{BV}}
\newcommand{\mub}{\overline{\mu}}
\newcommand{\muhat}{\hat{\mu}}

\newcommand{\hypref}[2]{\hyperref[#2]{#1 \ref*{#2}}}
\newcommand{\hypreff}[1]{\hyperref[#1]{(\ref*{#1})}}

\newcommand{\ob}[1]{^{(#1)}}

\newcommand{\xh}{\Xi}
\newcommand{\oh}[1]{O\left(#1\right)}
\newcommand{\xn}{\hat{x}}
\newcommand{\yn}{\hat{y}}
\newcommand{\On}{\hat{\Omega}}

\newcommand{\hn}{\hat{N}}

\begin{abstract}{}
Smooth solutions of the equation
\[
\op{div}\, \Bigg\{ \frac{g'\big(|\nabla u|\big)}{|\nabla u|}  \nabla u \Bigg\} = 0 
\]
are considered generating nonparametric $\mu$-surfaces in $\rz^3$, whenever $g$ is a function of linear growth satisfying in addition
\[
\int_0^\infty s g''(s) \D s < \infty \, .
\]
Particular examples are $\mu$-elliptic energy densities $g$ with exponent $\mu > 2$ (see \cite{BF:2003_1})
and the minimal surfaces belong to the class of $3$-surfaces.\\

Generalizing the minimal surface case we prove the closedness of a suitable differential form
$\hn \wedge \D X$. As a corollary we find an asymptotic conformal parametrization generated by this differential form.
\end{abstract}

\vspace*{2ex}

\centerline{\emph{Dedicated to S.~Hildebrandt}}

\vspace*{2ex}

\parindent0ex
\section{Introduction}\label{intro}
Starting with the variational point of view we like to mention three scenarios for variational problems with linear growth conditions.\\
 
The most prominent representative is the area minimizing problem
\begin{eqnarray}\label{intro 1}
J[\nabla w]  &:=& \iom F\big(\nabla w\big)\dx \to \min\, ,\nonumber \\[2ex]
F(Z) &:=&\sqrt{1+|Z|^2} \, , \quad Z\in \rz^2\, ,
\end{eqnarray}
within a suitable class of comparison functions $w$: $\rz^2 \supset \Omega \to \rz$. \\

The are uncountably many contributions to the study of problem \reff{intro 1}. 
We refer to the monographs \cite{Ni:1989_1}, \cite{Os:1986_1}, \cite{Gi:1984_1}, \cite{AFP:2000_1} and \cite{DHS:2010_1} 
giving a detailed picture of this classical problem.\\

There is another well known application involving variational problems with linear growth conditions: the theory
of perfect plasticity. We just mention \cite{Se:1985_1} as one of a series of papers written by Seregin and the
monograph of Temam, \cite{Te:2018_1} as well as the book \cite{FS:2000_1}.\\

As a third class of variational problems with linear growth conditions we like to mention  the discussion of
$\mu$-elliptic integrands with linear growth  introduced in  \cite{BF:2003_1}. 
Depending on the parameter $\mu$, this family includes the minimal surface case as one example with exponent $\mu=3$ 
and an approximation of  perfect plasticity is covered for large values of $\mu$.\\
 
In \cite{BF:2003_1} and in subsequent papers the question of existence and regularity of eventually relaxed solutions was studied  w.r.t.~to different 
circumstances. These investigations include also aspects from image analysis (see, e.g., \cite{BF:2012_1}), where the model serves as an 
appropriate TV-approximation. 
For an overview on the aspects of existence, relaxation and regularity of solutions under
the assumption of $\mu$-ellipticity we refer to \cite{Bi:1818} or \cite{BS:2013_1}.\\

While in the case of perfect plasticity and related applications the dual problem formulated in terms the stress tensor plays the key role, 
various geometric features are developed in the minimal surface case. Here the Euler-Lagrange equation for 
$C^2$-solutions of the variational problem \reff{intro 1}, i.e.~the minimal surface equation
\begin{equation}\label{intro 2}
u_{xx} \big(1+u_y^2\big) + u_{yy}\big(1+u_x^2\big) - 2 u_xu_y u_{xy} = 0 
\end{equation}
serves as a prototype for the study of elliptic PDEs arising in connection with problems in geometry.
In \cite{Si:1997_1} the reader will find an exposition with a particular focus to  
the geometric structure of equation \reff{intro 2}.\\

Our note on geometrical properties of what we call nonparametric $\mu$-surfaces is also strongly influenced by the pioneering work
on the minimal surface equation. To be precise, we consider the following theorem formulated by Dierkes, Hildebrandt and Sauvigny 
in the notion of differential forms.
\begin{theorem}\label{dhs} (\cite{DHS:2010_1}, Theorem Section 2.2)
A nonparametric surface $X(x,y) = \big(x,y,z(x,y)\big)$, described by a function $z=z(x,y)$ of class $C^2$ on a simply connected domain 
$\Omega$ of $\rz^2$, with the Gauss map $N=(\xi,\eta,\zeta)$ is a minimal surface if and only if the vector-valued differential form
$N \wedge \D X$ is a total differential, i.e., if and only there is a mapping $X^* \in C^2(\Omega,\rz^3)$ such that
\begin{equation}\label{intro 3}
-\D X^* = N \wedge \D X \, .
\end{equation}
If we write
\[
X^* = (a,b,c)\, , \quad N \wedge \D X = (\alpha , \beta ,\gamma)\, ,
\]
equation \reff{intro 3} is equivalent to
\[
-\D a = \alpha \, ,\quad -\D b = \beta \, ,\quad - \D c =\gamma \, .
\]
\end{theorem}

The particular importance of this theorem is evident by the fact,
that it serves as the main tool to prove that solutions of the minimal surface equation are analytic functions
and that $X^*$ induces a diffeomorphism leading to a conformal representation.\\

Motivated by the $\mu$-elliptic examples with linear growth mentioned above, we are faced with the question,
whether $C^2$-solutions of the corresponding Euler equations can also be characterized
by the closedness of suitable differential forms.

\begin{remark}\label{intro rem 1}
Let us shortly clarify the notion "nonparametic $\mu$-surface'': $\mu$-elliptic energy densities $g_\mu$ introduced in Example \ref{ex example 2} 
provide a typical motivation
for our studies. Going through the details of the proofs it becomes evident, that we just have to consider $C^2$-solutions of \reff{main 6}
together with our assumption \reff{main 5a}, which roughly speaking corresponds to the case $\mu >2$.\\

The limit case $\mu = 2$ also plays an important role in studying the regularity of solutions. While our geometric considerations are based on the finiteness
in condition \reff{main 5a}, we note that, e.g., in (1.9) of  \cite{BBM:2018_1}  the condition
\[
\int_1^\infty s g''(s) \D s = \infty
\]

characterizes the existence of regular solutions taking the boundary data of a Dirichlet problem. 
We also like to refer to the introductory remarks of \cite{BBM:2018_1} and to the classical paper \cite{NN:1959_1},
where related conditions can be found.

\end{remark}

Now let us introduce a more precise notation: given a simply connected domain $\Omega \subset \rz^2$ and 
a $C^2$-function $u$: $\Omega \to \rz$
we consider the nonparametric surface $X$: $\Omega \to \rz^3$
\[
X(x,y) = \big(x,\, y,\, u(x,y)\big)\, ,\quad (x,y) \in \Omega\, ,
\]

and denote the asymptotic normal by
\[
\hn  = (\hn_1 \, \hn_2\, \hn_3)\, ,\quad (x,y) \in \Omega\, ,
\]

with components
\begin{equation}\label{main 1}
\hn_1 = - \xh\big(|\nabla u|\big) u_x\, , \quad \hn_2 = - \xh\big(|\nabla u|\big) u_y\, ,
\quad \hn_3 = \xh\big(|\nabla u|\big) + \vartheta\big(|\nabla u|\big) \, .
\end{equation}
Here we let for $g$: $[0,\infty) \to \rz$  and  all $t \geq 0$
($g \in C^2\big([0,\infty)\big)$, $g'(0) = 0$, $g''(t) > 0$ for all
$t >0$) 
\begin{eqnarray}\label{main 2}
\xh(t)& :=&  \frac{g'(t)}{t}\, ,\\[2ex]
\label{main 3}
\vartheta(t)& :=& g(t) - tg'(t) - \xh(t) \, . 
\end{eqnarray}

The main hypothesis throughout this paper are summarized in the following Assumption.
\begin{ass}\label{main ass 1}
Let $g$: $[0,\infty) \to \rz$ be a function of class $C^2\big([0,\infty)\big)$ such that $g'(0) = 0$,
$g''(t) > 0$ for all $ t> 0$ and
such that with real numbers $a$, $A> 0$, $b$, $B\geq 0$
\begin{equation}\label{main 4}
a t - b \leq g(t) \leq A t + B \quad\mbox{for all}\quad t \geq 0 \, .
\end{equation}

Moreover, suppose that we have
\begin{equation}\label{main 5a}
\int_0^\infty s g''(s) \D s < \infty \, .
\end{equation}
\end{ass}

Before these hypotheses are discussed more detailed in Remark \ref{cor rem 1}, we state our main theorem:

\begin{theorem}\label{main}
With the notation introduced above we suppose that the general Assumption \ref{main ass 1} is valid, 
in particular we have \reff{main 5a}. Let $u$ denote a function of class $C^2(\Omega)$.\\

Then $u$ is a solution of 
\begin{equation}\label{main 6}
\op{div}\, \Bigg\{ \frac{g'\big(|\nabla u|\big)}{|\nabla u|}  \nabla u \Bigg\} = 0 
\end{equation}
on a simply connected domain $\Omega \subset \rz^2$ if and only if the vector-valued differential form
satisfies
\begin{equation}\label{main 7}
\hn \wedge \D X = -\D X^*
\end{equation}
for some mapping $X^* =: (a,b,c) \in C^2\big(\Omega;\rz^3\big)$, i.e.~$N \wedge \D X=:(\alpha,\beta,\gamma)$ is a total differential.
\end{theorem}

Let us recall the geometrical meaning of the differential form $N \wedge \D X$ in the minimal surface case by assuming that the surface is given
w.r.t.~isothermal parameters $(v,w)$. In this case we have
\[
N\wedge \D X = X_w \D v - X_v \D w 
\]
and the closedness of $N \wedge \D X$ corresponds to the minimal surface characterization $\Delta X = 0$ whenever
$X$ is given in conformal parameters. In fact it turns out that $X^*$ generates a diffeomorphism
leading from the nonparametric representation to conformal parameters and as a consequence to the analyticity of the solutions and
to Bernstein's theorem.\\

As it is also emphasized in \cite{DHS:2010_1}, this approach is of explicit geometrical nature related to the particular kind
of surfaces under consideration, which is in contrast to the abstract application of Lichtenstein's mapping theorem
to ensure the existence of a conformal representation. We refer to Hildebrandt's beautiful overview \cite{Hi:2011_1}
and the relation to Plateau's problem including quite recent results with von der Mosel \cite{HM:2005_1}, \cite{HM:2009_1}
and Sauvigny \cite{HS:2010_1}.\\

In our setting we are not in the minimal surface case having the analyticity of solutions, i.e.~we cannot 
argue with the help of suitably defined holomorphic functions. Thus we change the point of view in the sense that we are not mainly interested
in conformal representations. As an application of our main theorem we are rather interested in the question,
which kind of representation is generated by the mapping $X^*$ constructed in Theorem \ref{main}.
This means that we are looking for some kind of natural parametrisations for nonparametric $\mu$-surfaces.\\
   
It will turn out as a corollary of our main theorem, that both of the conformality relations are perturbed 
by the same function $\Theta$ and that asymptotically
the conformality relations are recovered.\\

Before giving a precise statement of this corollary, we like to include some additional remarks on
our assumptions.

\begin{remark}\label{cor rem 1}
\begin{enumerate}
\item Observe that we have an equivalent formulation of the main hypothesis \reff{main 5a}:
an integration by parts gives
\[
\int_0^t s g''(s) \D s = s \cdot g'(s)\Big|^t_0 - \int_0^t g'(s) \D s = t g'(t) - g(t) + g(0) \, .
\]

and we may write \reff{main 5a} in the form
\[
\lim_{t \to \infty} \Big[g(t) - t g'(t)\Big] = K:= g(0) -  \int_0^\infty s g''(s) \D s  \, .
\]

Replacing the function $g$ by the function $g-K$ (not changing equation \reff{main 6}), we may replace w.l.o.g.~assumption \reff{main 5a}
by 
\begin{equation}\label{main 5}
\lim_{t \to \infty} \Big[g(t) - t g'(t)\Big] = 0 
\end{equation}

and in the following \reff{main 5} can be taken as general assumption.

\item The convexity of $g$ immediately gives for all $t\geq 0$
\begin{equation}\label{cor 1}
g(0)  \geq g(t) - t g'(t) \, .
\end{equation}

Moreover, the function $g(t) - t g'(t)$ is a decreasing function in $[0, \infty)$ since we have for all $t \geq 0$
\begin{equation}\label{cor 2}
\frac{\D}{\D t} \Big[g(t) - t g'(t)\Big] = -t g''(t) \leq 0 \, .
\end{equation}

Thus,  by \reff{cor 1}, \reff{cor 2} and \reff{main 5}
\begin{equation}\label{cor 3}
g(0) \geq g(t) - t g'(t) \geq 0 \quad\mbox{for all}\quad t \geq 0 \, ,
\end{equation}

\item From the convexity of $g$ and the linear growth in the sense of \reff{main 4} we obtain the existence of
\[
\lim_{t\to \infty} g'(t) =: g_\infty' = \lim_{t\to \infty}\frac{g(t)}{t}\, .
\]

Writing 
\[
\frac{g(t) g'(t)}{t} - \big(g'(t)\big)^2 = R(t) \, , \quad\mbox{i.e}\quad
g'(t) \Big[g(t) - t g'(t)\Big] = t R(t) \, ,
\]

assumption \reff{main 5} in addition gives
\begin{equation}\label{cor 4}
R(t) = o\big(t^{-1}\big) \, .
\end{equation}

\item W.l.o.g.~let us suppose $g'_\infty =1$ and define for $t \geq 0$ the functions
\begin{equation}\label{cor 5}
g'(t) =: 1- h(t)\, ,\quad \xh(t) = \frac{1}{t} \big[1-h(t)\big] = \frac{g'(t)}{t} \geq 0  \, .
\end{equation}

Then $0\leq h(t) < 1$ for all $t\in \rz$ and $\lim_{t\to \infty} h(t) = 0$.
\end{enumerate}
\end{remark}

Now we formulate
\begin{cor}\label{corollary 1}
Suppose that $g$: $[0,\infty) \to \rz$ satisfies Assumption \ref{main ass 1} and suppose that
$u$: $\rz^2\supset \Omega \to \rz$ is a $C^2$-solution of \reff{main 6}, where $\Omega$ now in addition is assumed to
be convex.\\

Then the function $X^*$ described in Theorem \ref{main} generates an
asymptotically conformal parametrization $\chi$: $\On \to \rz^3$ of the surface $\op{graph} u$ in the following sense:

\begin{enumerate}
\item There is a function $\Theta$: $[0,\infty) \to \rz$ such that
\begin{eqnarray*}
\partial_{\xn} \chi \cdot \partial_{\yn} \chi &=& \frac{1}{\big(\op{det}D\Lambda\big)^2} \Theta\big(|\nabla u|\big) u_x u_y\, ,\\[2ex]
|\partial_{\xn}  \chi|^2 - |\partial_{\yn} \chi|^2 &=&  \frac{1}{\big(\op{det}D\Lambda\big)^2} \Theta\big(|\nabla u|\big) \big[u^2_x - u^2_y\big]\, .
\end{eqnarray*}

Here the diffeomorphism  $\Lambda$: $\Omega \to \On \subset \rz^2$ is given by
\begin{equation}\label{cor 6}
\Lambda(x,y) = \left(\begin{array}{c}x\\y\end{array}\right) +
 \left(\begin{array}{c}b(x,y)\\[1ex]-a(x,y)\end{array}\right)\, , \quad (x,y) \in \Omega \, ,
\end{equation}

and we define $\chi$: $\On \to \rz^3$ by
\begin{equation}\label{cor  7}
\chi: \, (\xn,\yn) \mapsto \Big(\Lambda^{-1}(\xn,\yn), u \circ \Lambda^{-1}(\xn,\yn) \Big)\, .
\end{equation}

\item There is a constant $c> 0$ such that for all $(x,y) \in \Omega$
\begin{equation}\label{cor 8}
{\rm det} \, D \Lambda \geq c \big(1+|\nabla u|\big)\, .
\end{equation}

\item If we suppose for some $\mu >2$ that for all $t \geq 0$
\begin{equation}\label{cor 9}
\big| g(t) - t g'(t)\big| \leq c_1 t^{2-\mu} \, ,\quad 
0 \leq 1-g'(t)  \leq  c_2 t^{1-\mu} \, ,
\end{equation}

with constants $c_1$, $c_2>0$, then we have for all $ t \gg 1$
\begin{equation}\label{cor 10}
|\Theta(t)| \leq  d_1 t^{2-\mu} + d_2 t^{-1}
\end{equation}
with some real numbers $d_1$, $d_2 >0$.
\end{enumerate}
\end{cor}

Let us close this introduction with three intuitive examples we have in mind.

\begin{example}\label{ex example 1}
The most prominent one is the minimal surface example given by
\[
g_{\min}(t) =  \sqrt{1+t^2} \, ,
\]
for which
\begin{eqnarray*}
\xh_{\min}(t) &=& \frac{1}{\sqrt{1+t^2}}\, , \\[2ex]
g_{\min}(t) - t g_{\min}'(t) &=& \frac{1}{\sqrt{1+t^2}} = \xh_{\min}(t)\, ,\\[2ex]
h_{\min}(t) &=&  1 - \frac{t}{\sqrt{1+t^2}} = \frac{1}{t\sqrt{1+t^2} + (1+t^2)}\, .
\end{eqnarray*}

With the help of this example we can always check our results by comparing to the classical ones.\\
\end{example}

\begin{example}\label{ex example 2}
We fix $\mu > 1$, $\mu\not= 2$,  and consider
\begin{equation}\label{ex 1}
g_\mu (t) = t + \frac{1}{\mu - 2} (1+t)^{2- \mu}  \, , \quad t \geq 0 \, .
\end{equation}

The particular choice $g_\mu$ is suitable for image analysis problems discussed in \cite{BF:2012_1}, since
for large $\mu$ we have convergence to the total variation energy.\\

For $g_\mu$ we have
\[
\xh_{\mu}(t) = \frac{1}{t} \Big[1- (1+t)^{1-\mu}\Big] \, , \quad
h_\mu (t) = (1+t)^{1-\mu}\, .
\]

We note that the growth of $h_\mu$ corresponds to the minimal surface case if $\mu =3$.
Moreover, we may use direct calculations to obtain
\begin{eqnarray*}
g_\mu (t) - t g_\mu'(t) &=& t + \frac{1}{\mu -2} (1+t)^{2-\mu} - t \Big[1-(1+t)^{1-\mu}\Big]\\[2ex]
&=& (1+t)^{1-\mu}\Bigg[\frac{1}{\mu-2} (1+t)+t\Bigg]\\[2ex]
&=& (1+t)^{1-\mu} \Bigg[\frac{\mu -1}{\mu -2} t + \frac{1}{\mu -2}\Bigg] \, .
\end{eqnarray*}

Hence, condition \reff{main 5} is satisfied whenever $\mu > 2$.\\
\end{example}

\begin{example}\label{ex example 3}
 In \cite{BF:2003_1} a variant of the above family is introduced by letting (again $\mu > 1$ is fixed)
 \[
 \hat{g}_\mu(t)  := \int_0^{t} \int_0^{s} \big(1+\tau^2\big)^{-\frac{\mu}{2}}\D \tau \D s \, , \quad t \geq 0 \, .
 \]
 
 This variant is of particular interest in our context since 
 the case $\mu =3$ exactly corresponds to the minimal surface case.\\
 
Here we immediately verify \reff{main 5a} if we again have  $\mu > 2$.
\end{example}

\section{The $\mu$-surface equation}\label{eq}

We suppose that a $F$: $\rz^2 \to \rz$ is given by
\[
F(Z) = g \big(|Z|\big) \, ,
\] 

where the function $g$: $[0,\infty) \to \rz$ satisfies the main Assumption \ref{main ass 1}. We observe that
\[
\nabla F(Z) =\frac{ g'\big(|Z|\big)}{|Z|} Z = \xh\big(|Z|\big) Z \, .
\]

and note that the variational problem
\[
J[w] := \iom F\big(\nabla w\big) \dx \to \min
\]
w.r.t.~a suitable class of comparison functions leads to the Euler equation
\begin{equation}\label{eq 1}
{\rm div}\, \Big\{ \nabla F( \nabla u) \Big\} = {\rm div}\, \Big\{ \xh\big( |\nabla u| \big) \nabla u \Big\} = 0 \, .
\end{equation}

Here and in the following we suppose that we have a solution $u$: $\rz^2\supset \Omega \to\rz$ 
to \reff{eq 1} which is at least of class $C^2(\Omega)$, $\Omega$ denoting an open set in $\rz^2$.\\

We write \reff{eq 1} in an explicit way:
\begin{eqnarray*}
\lefteqn{\Delta u \, \xh\big(|\nabla u|\big) + u_x \partial_x \xh\big(|\nabla u|\big) + u_y \partial_y \xh\big(|\nabla u|\big)}\\[2ex]
&=& \Delta u\,  \xh\big(|\nabla u|\big) \\[2ex]
&& + u_x \Bigg[\frac{\xh'\big(|\nabla u|\big)}{|\nabla u|} (u_x u_{xx} + u_y u_{xy})\Bigg]
+  u_y \Bigg[\frac{\xh'\big(|\nabla u|\big)}{|\nabla u|} (u_x u_{xy} + u_y u_{yy})\Bigg] \, ,
\end{eqnarray*}

which shows
\begin{eqnarray}\label{eq 2}
\op{div} \Big\{ \xh\big( |\nabla u| \big) \nabla u\Big\} &=&
u_{xx}\Bigg[\xh\big(|\nabla u|\big) + \frac{\xh' \big(|\nabla u|\big)}{|\nabla u|} u_x^2\Bigg]\nonumber\\[2ex]
&&+ u_{yy} \Bigg[ \xh\big(|\nabla u|\big) + \frac{\xh'\big(|\nabla u|\big)}{|\nabla u|} u_y^2 \Bigg]\nonumber\\[2ex]
&& + u_{xy}u_xu_y 2 \frac{\xh'\big(|\nabla u|\big)}{|\nabla u|}\, .
\end{eqnarray}

\vspace*{2ex}
\begin{example}\label{eq example 1}
In the minimal surface case we have the expressions
\begin{eqnarray*}
\xh_{\min} \big(|\nabla u|\big) + \frac{\xh'_{\min}\big(|\nabla u|\big)}{|\nabla u|} u_x^2
&=&  \frac{1+u_y^2}{\big(1+|\nabla u|^2\big)^{3/2}} \, ,\\[2ex]
\xh_{\min}\big(|\nabla u|\big) +  \frac{\xh'_{\min}\big(|\nabla u|\big)}{|\nabla u|} u_y^2
&=& \frac{1+u_x^2}{\big(1+|\nabla u|^2\big)^{3/2}} \, ,\\[2ex]
2 \frac{\xh'_{\min}\big(|\nabla u|\big)}{|\nabla u|}&=& -2 \frac{1}{\big(1+|\nabla u|^2\big)^{3/2}} \, ,
\end{eqnarray*}

thus \reff{eq 2} is equivalent to equation \reff{intro 2}.
\end{example}

\section{Proof of the main theorem}\label{Potential}
With the definition of the asymptotic normal $\hn$ one computes  $\hn \wedge \D X = (\alpha ,\beta ,\gamma)$ with components
\begin{eqnarray}\label{pm 1}
\alpha &=& - \xh\big(|\nabla u|\big)  u_y \du - \Big[\xh\big(|\nabla u|\big) +\vartheta\big(|\nabla u|\big) \Big]\D y\nonumber\\[2ex]
\beta &=& \Big[\xh\big(|\nabla u|\big) + \vartheta\big(|\nabla u|\big)\Big]\D x + \xh\big(|\nabla u|\big) u_x\D u\nonumber\\[2ex]
\gamma &=& \xh\big(|\nabla u|\big) u_y \D x - \xh\big(|\nabla u|\big) u_x \D y\, .
\end{eqnarray}

We observe
\[
\D u = u_x \D x + u_y \D y
\]

and obtain
\begin{eqnarray}\label{pm 2}
\alpha &=& - \xh\big(|\nabla u|\big)  u_y u_x\dx - \xh\big(|\nabla u|\big)u_y u_y\D y
- \Big[\xh\big(|\nabla u|\big) +\vartheta\big(|\nabla u|\big) \Big]\D y\nonumber\\[2ex]
&=&  - \xh\big(|\nabla u|\big) u_x u_y\D x - \Big[\xh\big(|\nabla u|\big) \big(1+u_y^2\big) + \vartheta\big(|\nabla u|\big)\Big]\D y \, .
\end{eqnarray}

In the same way we get 
\begin{equation}\label{pm 3}
\beta =  \Big[\xh\big(|\nabla u|\big) \big(1+u_x^2\big) + \vartheta\big(|\nabla u|\big) \Big]\D x + \xh\big(|\nabla u|\big) u_x u_y\D y \, .
\end{equation}

From \reff{pm 1} -- \reff{pm 3} we obtain the equation
\begin{eqnarray}\label{pm 4}
\D \alpha &=& \partial_y \Big[ \xh\big(|\nabla u|\big) u_x u_y\Big] \D x\wedge \D y \nonumber\\[2ex]
&&- \partial _x\Big[ \xh\big(|\nabla u|\big) \big(1+u_y^2\big)+ \vartheta \big(|\nabla u|\big)\Big] \D x\wedge \D y\nonumber\\[2ex]
&=:& \big[\partial_y \psi_1 -  \partial_x \psi_2\big] \D x \wedge\D y\, .
\end{eqnarray}

The exterior derivative of the form $\beta$ is given by 
\begin{eqnarray}\label{pm 5}
\D \beta &=& - \partial _y\Big[ \xh\big(|\nabla u|\big) \big(1+u_x^2\big)+ \vartheta \big(|\nabla u|\big)\Big] \D x\wedge \D y\nonumber\\[2ex]
&&+ \partial_x \Big[ \xh\big(|\nabla u|\big) u_x u_y\Big] \D x \wedge\D y \nonumber\\[2ex]
&=:&\big[- \partial_y \varphi_1 + \partial_x \varphi_2\Big] \D x\wedge\D y \, .
\end{eqnarray}

Finally we have for $\D \gamma$
\begin{equation}\label{pm 6}
\D \gamma = - \partial_y \Big[\xh\big(|\nabla u|\big) u_y\Big] \D x \wedge \D y -
\partial_x \Big[\xh\big(|\nabla u|\big) u_x\Big] \D x \wedge \D y\, .
\end{equation}

Let us first consider \reff{pm 6} by computing
\begin{eqnarray*}
\lefteqn{\partial_y \Big[ \xh\big(|\nabla u|\big)u_y \Big] 
+ \partial_x \Big[ \xh\big(|\nabla u|\big)u_x \Big]}\\[2ex] 
&=& u_{yy} \xh\big(|\nabla u|\big) + u_y \xh'\big(|\nabla u|\big) \frac{u_x u_{xy}+u_y u_{yy}}{|\nabla u|}\\[2ex]
&& +  u_{xx} \xh\big(|\nabla u|\big) + u_x \xh'\big(|\nabla u|\big) \frac{u_x u_{xx}+u_y u_{xy}}{|\nabla u|}\\[2ex]
&=& u_{xx} \Bigg[\xh\big(|\nabla u|\big) + \frac{\xh'\big(|\nabla u|\big)}{|\nabla u|} u_x^2\Bigg]\nonumber+ u_{yy} \Bigg[\xh\big(|\nabla u|\big) + \frac{\xh'\big(|\nabla u|\big)}{|\nabla u|} u_y^2\Bigg]\\[2ex]
&&+ u_{xy} \Bigg[2 u_xu_y \frac{\xh'\big(|\nabla u|\big)}{|\nabla u|}\Bigg] \, ,
\end{eqnarray*}

hence we have recalling \reff{eq 2}
\begin{equation}\label{pm 7}
\D \gamma = 0 \, .
\end{equation}

Next we discuss \reff{pm 5}. Direct calculations show
\begin{eqnarray*}
\partial_y \varphi_1 & = & \frac{\xh'\big(|\nabla u|\big)}{|\nabla u|} (u_xu_{xy} + u_y u_{yy}) \big(1 + u_x^2\big)+ 
\xh\big(|\nabla u|\big) 2u_x u_{xy}\\[2ex]
&& + \frac{\vartheta' \big(|\nabla u|\big)}{|\nabla u|}  (u_x u_{xy}+ u_{y} u_{yy})
\end{eqnarray*}

and in addition
\[
- \partial_x \varphi_2 = - \frac{\xh' \big(|\nabla u|\big)}{|\nabla u|} (u_x u_{xx} + u_{y}u_{xy}) u_x u_y
- \xh\big( |\nabla u|\big) (u_{xx} u_y + u_x u_{xy}) \, .
\]

Combining these equations one obtains
\begin{eqnarray}\label{pm 8}
\partial_y \varphi_1 - \partial_x \varphi_2 &=& - u_y u_{xx}\Bigg[\xh\big(|\nabla u|\big)+ \frac{\xh'\big(|\nabla u|\big)}{|\nabla u|} u_x^2\Bigg]
\nonumber\\[2ex]
&& -u_y u_{yy} \Bigg[- \frac{\xh'\big(|\nabla u|\big)}{|\nabla u|}\big(1+u_x^2\big) - \vartheta'\big( |\nabla u|\big) \Bigg]\nonumber\\[2ex]
&& - u_{xy} \Bigg[ - \frac{\xh'\big(|\nabla u|\big)}{|\nabla u|} u_x \big(1+u_x^2\big) - \xh\big(|\nabla u|\big) u_x \nonumber\\[2ex]
&& - \frac{\vartheta' \big(|\nabla u|\big)}{|\nabla u|} u_x
+u_x u_y^2 \frac{\xh'\big(|\nabla u|\big)}{|\nabla u|}\Bigg]\nonumber\\[2ex]
&=:& - u_y u_{xx} T_1 - u_y u_{yy} T_2 - u_{xy} T_3 \, .  
\end{eqnarray}

Now we compare \reff{pm 8} with the equation \reff{eq 2} and observe that $T_1$ is the first coefficient 
on the right-hand side of \reff{eq 2}. \\

If we can show in additon that 
\begin{eqnarray}\label{pm 9}
T_2 &=& \xh \big( |\nabla u|\big) + \frac{\xh'\big(|\nabla u|\big)}{|\nabla u|} u_y^2\, ,\\[2ex]
\label{pm 10}
T_3 &=& u_y \Bigg[ 2 u_x u_y  \frac{\xh'\big( |\nabla u|\big)}{|\nabla u|}\Bigg]\, .  
\end{eqnarray}

then we obtain
\begin{equation}\label{pm 11}
\partial_y \varphi_1 - \partial_x \varphi_2 = 0\, , \quad\mbox{hence}\quad \D \beta = 0 \, . 
\end{equation}

For discussing \reff{pm 9}, i.e.~the validity of the equation
\[
-\frac{\xh'\big(|\nabla u|\big)}{|\nabla u|} \big(1+u_x^2\big) - \frac{\vartheta'\big(|\nabla u|\big)}{|\nabla u|} = 
 \xh \big( |\nabla u|\big) + \frac{\xh'\big(|\nabla u|\big)}{|\nabla u|} u_y^2 \, ,
\]

we observe that the latter identity is equivalent to
\begin{equation}\label{pm 12}
- \frac{\vartheta' \big(|\nabla u|\big)}{|\nabla u|} = \xh \big(|\nabla u|\big) + \frac{\xh'\big(|\nabla u|\big)}{|\nabla u|}\big(1+ |\nabla u|^2\big) \, .
\end{equation}

The definition of $\vartheta$ (recall \reff{main 3}) now yields
\begin{eqnarray}\label{pm 13}
\vartheta'(t) &=& -  \Big[t \xh(t) + \xh'(t) \big(1+t^2\big)\Big] \, .
\end{eqnarray}

which immediately gives \reff{pm 12}.\\

Equation \reff{pm 10} takes the form
\begin{eqnarray*}
2 u_y^2 u_x \frac{\xh'\big(|\nabla u|\big)}{|\nabla u|} &=&
 - \frac{\xh'\big(|\nabla u|\big)}{|\nabla u|} u_x \big(1+u_x^2\big) - \xh\big(|\nabla u|\big) u_x \\[2ex]
&& - \frac{\vartheta' \big(|\nabla u|\big)}{|\nabla u|} u_x
+ u_x u_y^2 \frac{\xh'\big(|\nabla u|\big)}{|\nabla u|} \, ,
\end{eqnarray*}

i.e. 
\[
u_y^2 u_x \frac{\xh'\big(|\nabla u|\big)}{|\nabla u|} = 
 - \frac{\xh'\big(|\nabla u|\big)}{|\nabla u|} u_x \big(1+u_x^2\big) - \xh\big(|\nabla u|\big) u_x 
- \frac{\vartheta' \big(|\nabla u|\big)}{|\nabla u|} u_x 
\]

and this relation is equivalent to 
\[
u_x\Bigg[\frac{\xh'\big(|\nabla u|\big)}{|\nabla u|}\big(1+|\nabla u|^2\big) + \xh\big( |\nabla u|^2\big)\Bigg] = 
- \frac{\vartheta' \big(|\nabla u|\big)}{|\nabla u|} u_x \, .
\]

Again we end up with \reff{pm 13}, hence we also have \reff{pm 10}
and finally \reff{pm 11}.\\

It remains to consider \reff{pm 4} and to show
\begin{equation}\label{pm 14}
\D \alpha = 0 \, .
\end{equation}

However this can be done using the same arguments leading to \reff{pm 11}.\\

With \reff{pm 7}, \reff{pm 11} and \reff{pm 14} it is shown that $\hn \wedge \D X$ is a closed differential form
and we have \reff{main 7} on the simply connected domain $\Omega$, thus \mbox{Theorem \ref{main}} is established.\qed\\

\section{A parametrization generated by $X^*$}\label{param}
In the minimal surface case, Theorem \ref{main} implies the conformal representation and the analyticity of 
nonparametric minimal surfaces as outlined, e.g., in Section 2.3 of \cite{DHS:2010_1}. The approach given there
can be done by just varying the second surface-parameter.\\ 

Here we do not expect analytic solutions in general and we prefer to follow a variant given in, e.g., \cite{Os:1986_1}, which
yields Lemma \ref{conform lem 1} in its symmetric formulation.\\

Given $b(x,y)$, $a(x,y)$ according to Theroem \ref{main} we consider the differential form 
\[
\omega := b(x,y)\D x - a(x,y) \D y\, .
\]

Then we have on account of $-\D a = \alpha$, $-\D b = \beta$ and recalling \reff{pm 2}, \reff{pm 3}
\[
\D\omega = - b_y(x,y)\D x \wedge \D y - a_x(x,y)\D x\wedge \D y = 0\, ,
\]

which means that the form is closed and we may define the line integral
\begin{equation}\label{param 1}
E(x,y) := \int_{(x_0,y_0)}^{(x,y)}  \omega \, ,
\end{equation}

where we have
\begin{equation}\label{param 2}
\nabla E(x,y) = \left( \begin{array}{c}
b(x,y)\\[1ex]-a(x,y)
\end{array}\right) \, .
\end{equation}

With the notation of \reff{pm 4} and \reff{pm 5} one obtains
\begin{equation}\label{param 3}
\left(\begin{array}{c} \varphi_1\\[1ex] \varphi_2\end{array}\right) = \nabla b \, , \quad
\left(\begin{array}{c} \psi_1\\[1ex] \psi_2\end{array}\right) = - \nabla a \, .
\end{equation}

Combining \reff{param 2} and \reff{param 3} finally gives
\begin{equation}\label{param 4}
D^2 E = \left(\begin{array}{cc}
\varphi_1&\varphi_2\\[1ex]
\psi_1 & \psi_2 
\end{array}\right) \, .
\end{equation}

Discussing $D^2E$ the relevance of condition \reff{main 5} becomes obvious.

\begin{proposition}\label{param prop 1}
Consider the function $E$: $\Omega \to \rz$ given in \reff{param 1} and suppose that \reff{main 5} holds.\\

Then for any $(x,y) \in \Omega$ the bilinear form  $D^2 E(x,y)$ is positive definite.\\
\end{proposition}

\emph{Proof of Propositon \ref{param prop 1}.} Fix some $(x,y) \in \Omega$ and abbreviate $D^2E(x,y)$ through $D^2E$.
Observe that \reff{param 4} implies for all $\eta = (\eta_1,\eta_2) \in \rz^2-\{0\}$:
\begin{eqnarray}\label{param 5}
D^2E (\eta,\eta) &=& \xh\big(|\nabla u|\big) \Big[\eta_1^2 \big(1+u_x^2\big) + 2\eta_1\eta_2 u_x u_y
+ \eta_2^2 \big(1+u_y^2\big)\Big] \nonumber\\[2ex]
&& + \vartheta \big(|\nabla u|\big) \big(\eta_1^2+\eta_2^2\big)\, .
\end{eqnarray}

Considering \reff{param 5} we note 
\[
|2 \eta_1 \eta_2 u_x u_y| \leq \eta_1^2 u_x^2 + \eta_2^2 u_y^2 \, ,
\]

hence the definition \reff{main 3} of the function $\vartheta$ shows
\begin{eqnarray*}
D^2E (\eta,\eta) & \geq & \Big[\xh\big(|\nabla u|\big) + \vartheta \big(|\nabla u|\big)\Big] |\eta|^2\\[2ex]
&=& \Bigg[ g(|\nabla u|) - |\nabla u| g'\big(|\nabla u|\big)  \Bigg] |\eta|^2  \, .
\end{eqnarray*}

Thus we can apply hypothesis \reff{main 5} to see that $D^2 E$ is positive definite. \qed\\

In the next step we introduce a diffeomorphism generated by the gradient field $\nabla E$, i.e.~by $X^*$ (recall \reff{cor 6}) :
\begin{equation}\label{conform 1}
\Lambda(x,y) = \left(\begin{array}{c}x\\y\end{array}\right) + \nabla E(x,y) \, , \quad (x,y) \in \Omega \, .
\end{equation}

As shortly outlined in the appendix, well known arguments show that $\Lambda$ in fact is a diffeomorphism onto its image.\\

A more refined analysis of $\Lambda$ of course depends on the underlying function $X^*$. In the classical minimal
surface case we are directly lead to a conformal parametrization without referring to Lichtenstein's mapping theorem
(see, e.g., \cite{DHS:2010_1}, Section 2.3 for further comments).\\

Here we expect the asymptotic correspondence to this method as a natural consequence of our main theorem \ref{main}.\\ 

We start with an estimate for the Jacobian which proves the claim \reff{cor 8} of Corollary \ref{corollary 1}.
\begin{proposition}\label{conform prop 1}
Suppose that we have Assumption \ref{main ass 1} and consider the diffeomorhism $\Lambda$ defined in \reff{conform 1}.\\

Then we have for all $(x,y)\in \Omega$
\[
{\rm det} \, D \Lambda \geq 1 +\xh \big(|\nabla u|\big) \big(1+|\nabla u|^2\big)\, .
\]
\end{proposition}

\emph{Proof of Proposition \ref{conform prop 1}.} The Jacobian of $\Lambda$ is given by
\begin{eqnarray}\label{conform 2}
\op{det} D \Lambda &=& \partial_x \Lambda_1 \partial_y \Lambda_2 - \partial_x \Lambda_2\partial_y \Lambda_1\nonumber\\[2ex]
&=& \Big[1+ \big[\xh\big(|\nabla u|\big)(1+u_x^2) + \vartheta\big(|\nabla u|\big)\big]\Big]\nonumber\\[2ex]
&& \quad \cdot  \Big[1+ \big[\xh\big(|\nabla u|\big)(1+u_y^2) + \vartheta\big(|\nabla u|\big)\big]\Big]\nonumber\\[2ex]
&&- \xh^2\big(|\nabla u|\big) \big(|\nabla u|\big)u^2_xu^2_y\nonumber\\[2ex]
&=& 1+ \xh\big(|\nabla u|\big)  \big(2+|\nabla u|^2\big) + \xh^2\big(|\nabla u|\big) \big(1+|\nabla u|^2\big)\nonumber\\[2ex]
&&+ \vartheta\big(|\nabla u|\big) + \vartheta\big(|\nabla u|\big)\xh\big(|\nabla u|\big) \big(2+|\nabla u|^2\big)\nonumber\\[2ex]
&&+ \vartheta^2\big(|\nabla u|\big) \, .
\end{eqnarray}

Moreover, the definition \reff{main 3} of $\vartheta$  shows
\begin{eqnarray}\label{conform 3}
\lefteqn{\vartheta\big(|\nabla u|\big) + \vartheta\big(|\nabla u|\big)\xh\big(|\nabla u|\big) \big(2+|\nabla u|^2\big)
+ \vartheta^2\big(|\nabla u|\big)}\nonumber\\[2ex]
&=  & \Big[g\big(|\nabla u|\big) - |\nabla u| g'\big(|\nabla u|\big)\Big] - \xh\big(|\nabla u|\big) \nonumber\\[2ex]
&&+ \Big[g\big(|\nabla u|\big) - |\nabla u| g'\big(|\nabla u|\big)\Big]\xh\big(|\nabla u|\big)\big(2+|\nabla u|^2\big) \nonumber\\[2ex]
&& - \xh^2\big(|\nabla u|\big)\big(2+|\nabla u|^2\big) + \Big[ g\big(|\nabla u|\big) - |\nabla u| g'\big(|\nabla u|\big)\Big]^2\nonumber\\[2ex]
&& - 2 \xh\big(|\nabla u|\big) \Big[g\big(|\nabla u|\big) - |\nabla u|g'\big(|\nabla u|\big) \Big]+ \xh^2\big(|\nabla u|\big) \, .
\end{eqnarray}

Combining \reff{conform 2} and \reff{conform 3} yields
\begin{eqnarray}\label{conform 4}
\op{det} D \Lambda &=& 1+ \xh\big(|\nabla u|\big)  \big(1+|\nabla u|^2\big) \nonumber\\[2ex]
&&+ \Big[g\big(|\nabla u|\big) - |\nabla u| g'\big(|\nabla u|\big)\Big]\nonumber\\[2ex]
&& \cdot \Bigg[1+ \xh\big(|\nabla u|\big) |\nabla u|^2 +  \Big[g\big(|\nabla u|\big) - |\nabla u| g'\big(|\nabla u|\big)\Big]\Bigg]\, .
\end{eqnarray}

With \reff{conform 4} the proof of Proposition \ref{conform prop 1} again follows from \reff{main 5}. \qed\\

Let us now rewrite the functions $\varphi_1$ and $\psi_2$ in the following form
\begin{eqnarray}\label{conform 5}
\varphi_1 &=& \xh\big(|\nabla u|\big)\big(1+u_x^2\big) + g\big(|\nabla u|\big) - \xh\big(|\nabla u|\big) |\nabla u|^2
-\xh\big(|\nabla u|\big)\nonumber\\[2ex]
&=& g\big(|\nabla u|\big) - \xh\big(|\nabla u|\big) u_y^2 \, ,\nonumber\\[2ex]
\psi_2 &=&  g\big(|\nabla u|\big) - \xh\big(|\nabla u|\big) u_x^2 \, .
\end{eqnarray}

Using \reff{conform 5} we recall the definition of $\Lambda$ and \reff{param 4}, hence
\[
D \Lambda =  
\left(\begin{array}{cc}
1+ \varphi_1 & \varphi_2\\[4ex]
\psi_1 & 1+ \psi_2 
\end{array}\right)\, ,
\]

and we calculate for all $(\xn , \yn)\in \On$ ($(x,y) := \lambda^{-1}(\xn,\yn)$)
\begin{eqnarray}\label{conform 6}
\lefteqn{D\big(\Lambda^{-1}\big)(\xn,\yn) = \big(D\Lambda(x,y)\big)^{-1}}\nonumber\\[4ex]
&=& \left(\begin{array}{cc}
1+ g \big(|\nabla u|\big) -\xh\big(|\nabla u|\big)u_y^2& \xh\big(|\nabla u|\big) u_x u_y\\[4ex]
\xh\big(|\nabla u|\big) u_x u_y & 1+ g\big(|\nabla u|\big) -\xh\big(|\nabla u|\big)u_x^2
\end{array}\right)^{-1} \nonumber\\[4ex]
&=&\frac{1}{\op{det} D\Lambda} \left(\begin{array}{cc}
1+ g\big(|\nabla u|\big) -\xh\big(|\nabla u|\big)u_x^2& -\xh\big(|\nabla u|\big) u_x u_y\nonumber\\[4ex]
-\xh\big(|\nabla u|\big) u_x u_y & 1+ g\big(|\nabla u|\big) -\xh\big(|\nabla u|\big)u_y^2
\end{array}\right)\nonumber \\[4ex]
&=:& \frac{1}{\op{det} D\Lambda}\,  \Pi(x,y) \, .
\end{eqnarray}

In addition we compute
\begin{eqnarray}\label{conform 7}
D\Big(u \circ \Lambda^{-1}\Big)(\xn ,\yn) &=& Du(x,y) \, D\big(\Lambda^{-1}\big)(\xn,\yn)\nonumber\\[2ex]
&=& \frac{1}{\op{det}D\Lambda} \Pi(x,y) \, \nabla u(x,y)\nonumber\\[2ex]
&=& \frac{1}{\op{det}D\Lambda}\, \big[1+g - |\nabla u|^2 \xh\big] \nabla u\big]\nonumber\\[2ex]
&=:&\frac{1}{\op{det}D\Lambda}\, \pi(x,y)  \, .
\end{eqnarray}

After these preparations we consider the parametrization of the surface $\op{graph}(u)$ already introduced in \reff{cor 7}: 
\[
\chi: \, (\xn,\yn) \mapsto \Big(\Lambda^{-1}(\xn,\yn), u \circ \Lambda^{-1}(\xn,\yn) \Big)\, , \quad (\xn ,\yn) \in \On\, .
\]

Then \reff{conform 6} and \reff{conform 7} yield
\begin{equation}\label{conform 8}
D\chi (\xn,\yn) =: \frac{1}{\op{det} D\Lambda} \Big(X\,\, Y\big) = \frac{1}{\op{det} D\Lambda}
\left(\begin{array}{c}
\Pi (x,y)\\[1ex]
\pi (x,y) \end{array}\right) \, ,
\end{equation}

where \reff{conform 6}, \reff{conform 7} and the definition \reff{conform 8} of $X$, $Y$ imply
\begin{eqnarray*}
X = \left(\begin{array}{c}
1+g\big(|\nabla u|\big) - \xh\big(|\nabla u|\big) u_x^2\\[2ex] 
-\xh\big(|\nabla u|\big) u_x u_y\\[2ex] 
u_x\big[1+g\big(|\nabla u|\big) -\xh\big(|\nabla u|\big) |\nabla u|^2\big]
\end{array}\right) \, , \\[4ex]
Y = \left(\begin{array}{c}
-\xh\big(|\nabla u|\big) u_x u_y\\[2ex] 
1+g\big(|\nabla u|\big) - \xh\big(|\nabla u|\big) u_y^2\\[2ex] 
u_y\big[1+g\big(|\nabla u|\big) -\xh\big(|\nabla u|\big) |\nabla u|^2\big]
\end{array}\right) \, .
\end{eqnarray*}

Now we come to the last part of Corollary \ref{corollary 1}:
\begin{lemma}\label{conform lem 1}
If $X$ and $Y$ are given as above, then we have the equations
\begin{eqnarray*}
X \cdot Y &=& u_x u_y \Theta\big(|\nabla u|\big) \, ,\\[2ex]
|X|^2 - |Y|^2 &=& \big[u_x^2 -u_y^2\big] \Theta\big(|\nabla u|\big) \, ,
\end{eqnarray*}
where the function $\Theta$: $[0,\infty) \to \rz$ is given by
\[
\Theta(t) =
\Bigg[1 - \frac{g(t)g'(t)}{t}\Bigg]\nonumber +\Bigg[\big(g(t)-t g'(t)\big)-\frac{g'(t)}{t}\Bigg] \Big[2+ \big(g(t)-t g'(t)\big)\Big]\, .
\]
\end{lemma}

\emph{Proof of Lemma \ref{conform lem 1}.} By elementary calculations we obtain
\begin{eqnarray}\label{conform 9}
X \cdot Y &=& u_x u_y\Bigg[ - \xh\big(|\nabla u|\big)\Big[2\big(1+g\big(|\nabla u|\big) - \xh\big(|\nabla u|\big) |\nabla u|^2\Big]\nonumber\\[2ex]
&& + \Big[1+g\big(|\nabla u|\big) -\xh\big(|\nabla u|\big) |\nabla u|^2\Big]^2\Bigg]\nonumber\\[2ex]
&=:& u_x u_y \, \tilde{\Theta}\big(|\nabla u|\big) \, .
\end{eqnarray}

as well as
\begin{eqnarray}\label{conform 10}
|X|^2 - |Y|^2 &=& \big[u_x^2-u_y^2\big] \, \tilde{\Theta} \big(|\nabla u|\big)
\end{eqnarray}
with the same function $\tilde{\Theta}$.\\

Let us write the function $\tilde{\Theta} (t)$ from \reff{conform 9} and \reff{conform 10} as a more convenient form. We have
\[
2 \big(1+g(t)\big) -\xh(t) t^2 = \big(2+g(t)\big) + \big(g(t) - tg'(t)\big)\, ,
\]

hence
\begin{eqnarray}\label{conform 11}
\lefteqn{- \xh(t) \Bigg[2 \big(1+g(t)\big) -\xh(t) t^2\Bigg]}\nonumber\\[2ex]
& = &
- \frac{\big(2+g(t)\big)g'(t)}{t} - \frac{g'(t)}{t} \big(g(t)-t g'(t)\big) \, .
\end{eqnarray}

Moreover, we note
\begin{equation}\label{conform 12}
\Big[1+g(t)- t g'(t)\Big]^2 = 1 + 2 \big(g(t)-tg'(t)\big) + \big(g(t)-t g'(t)\big)^2 \, .
\end{equation}

Adding \reff{conform 11} and \reff{conform 12} we obtain
\begin{eqnarray*}
\tilde{\Theta}(t) &=&  \Bigg[1 - \frac{g(t)g'(t)}{t}\Bigg] - \frac{2 g'(t)}{t} + \Bigg[2-\frac{g'(t)}{t}\Bigg]\big(g(t)-t g'(t)\big)\nonumber\\[2ex]
&& + \big(g(t) - t g'(t)\big)^2\nonumber\\[2ex]
&=&
\Bigg[1 - \frac{g(t)g'(t)}{t}\Bigg]\nonumber\\[2ex] 
&&+\Bigg[\big(g(t)-t g'(t)\big)-\frac{g'(t)}{t}\Bigg] \Big[2+ \big(g(t)-t g'(t)\big)\Big] = \Theta(t)\, ,
\end{eqnarray*}

hence we have proved Lemma \ref{conform lem 1}. \qed\\

\begin{example}\label{conform example 2}
We recall that in the minimal surface case we have for all $t\geq 0$
\[
g(t) - tg'(t) = \frac{g'(t)}{t}\quad\mbox{and}\quad \frac{g(t) g'(t)}{t} \equiv 1 \, .
\]
\end{example}

As a consequence of Lemma \ref{conform lem 1} we obtain an explicit asymptotic expansion for the function $\Theta$.
\begin{cor}\label{conform cor 1} 
Suppose that we have the assumptions of Corollary \ref{corollary 1}. 
Then we have \reff{cor 10}, i.e.~for all $t > 0$
\[
|\Theta(t)| \leq  d_1 |t|^{2-\mu} + d_2|t|^{-1}
\]
with some real numbers $d_1$, $d_2 >0$.
\end{cor}

\emph{Proof of Corollary \ref{conform cor 1}.} As in Remark \ref{cor rem 1}, $ii$), we now have using \reff{cor 9}
\[
g'(t) \Big[g(t) - t g'(t)\Big] = t R(t) \, , \quad\mbox{i.e.}\quad R(t) = O\big(t^{1-\mu}\big) \, .
\]

By Remark \ref{cor rem 1}, $ii$), we obtain in addition
\begin{eqnarray*}
1- \frac{g(t) g'(t)}{t} &=& \Big[1-\big(g'(t)\big)^2\Big] + O\big(t^{1-\mu}\big)\\[2ex]
&=& h(t) \big(1+g'(t)\big) + O\big(t^{1-\mu}\big) \, .
\end{eqnarray*}

This, together with the boundedness of $g'$ and once more applying \reff{cor 9} gives the corollary. \qed\\

We finally note that Proposition \ref{conform prop 1}, Lemma \ref{conform lem 1} and Corollay \ref{conform cor 1}
together yield the proof of Corollary \ref{corollary 1}, \hfill\qed \\

\section{Appendix}\label{app}
For the sake of completeness we append Proposition \ref{app prop} which guarantees
the fact that the function $\Lambda$ defined in \reff{conform 1} is a diffeomorhism onto its image. 
The relation \reff{app 1} is of particular interest if $\Omega = \rz^2$. The short proof is outlined below following
Osserman's book \cite{Os:1986_1}.

\begin{proposition}\label{app prop}
With the hypotheses of Proposition \ref{param prop 1} we suppose that $\Omega$ is a convex open set and let 
$\Lambda$:  $\Omega \to G:= \Lambda(\Omega)$,
\[
\Lambda(x,y) = \left(\begin{array}{c}x\\y\end{array}\right) + \nabla E(x,y) \, , \quad (x,y) \in \Omega \, ,
\]
with $E$ defined in \reff{param 1}.\\

Then $\Lambda$ is a diffeomorphism, moreover we have for any $B_r(x_0,y_0) \subset \Omega$
\begin{equation}\label{app 1}
B_r \Big(\Lambda\big( (x_0 ,y_0)\big)\Big) \subset G \, .
\end{equation}
\end{proposition}

\emph{Proof of Proposition \ref{app prop}.}
On the convex set $B_r(0)$ and we let
\[
e(t) := E\big(tx+ (1-t)y\big)\, ,\quad t\in[0,1]\, .
\]

Then
\[
e''(t) = D^2E\big(tx+(1-t)y\big) \big(x-y,x-y\big) > 0 \quad \mbox{for all $t\in (0,1)$} \, ,
\]

which shows that $e'(t)$ is an increasing function, i.e.
\begin{equation}\label{app 2}
0 < e'(1) - e'(0) = \Big[\nabla E(x) - \nabla E(y)\Big]\cdot (x-y) \, . 
\end{equation}

By the definition of $\Lambda$ and by \reff{app 1}, \reff{app 2} we obtain ($x\not= y \in B_r(0)$)
\begin{eqnarray*}
\lefteqn{\big[\Lambda(x_1,y_1) - \Lambda(x_2.y_2)\big]\cdot \big[(x_1,y_1)-(x_2,x_2)\big]}\nonumber\\[2ex]
&=&| (x_1,y_1) -(x_2,y_2)|^2\nonumber\\[2ex]
&&+ \Big[\nabla E(x_1,y_1) - \nabla E(x_2,y_2)\Big]
\cdot \big[(x_1,y_1)-(x_2,x_2)\big] \\[2ex]
& >&| (x_1,y_1) -(x_2,y_2)|^2  \, ,
\end{eqnarray*}

and by the Cauchy-Schwarz inequality this implies
\begin{equation}\label{app 3}
\big| \Lambda(x_1,y_1) -\Lambda(x_2,y_2)\big| >| (x_1,y_1) -(x_2,y_2)|  \, .
\end{equation}

Observe that \reff{app 3} gives the injectivity of $\Lambda$ and since in addition 
\[
\op{det} D \Lambda (x,y) >0 \quad\mbox{for}\quad (x,y) \in B_r(0)\, , 
\]
$\Lambda$: $B_r(0) \to G := \Lambda\big(B_r(0)\big)$ is a global diffeomorphism and we have \reff{app 1}
if $G = \rz^2$.\\

In the case $G \not= \rz^2$ there exists $(\xn,\yn) \in \partial \big(\rz^2 \setminus G\big) = \partial G$, such that
 \[
\Big|(\xn,\yn) - \Lambda\big((x_0,y_0)\big)\Big| = \inf_{\eta \in \rz^2\setminus G} \Big|\eta - \Lambda\big((x_0,y_0)\big)\Big| \, .
 \] 

This means that we find sequences $\{(x,y)^{(k)}\}$, $\{(\xn ,\yn)^{(k)}=\Lambda\big((x,y)^{(k)}\big)\}$ such that as $k \to \infty$
\[
(x,y)^{(k)} \to (x,y) \in \partial B_r(x_0,y_0) \quad\mbox{and}\quad (\xn ,\yn)^{(k)} \to (\xn ,\yn) \, .
\]

The definition of $(\xn ,\yn)$ finally yields 
\[
\inf_{\eta\in \rz^2\setminus G} \Big|\eta - \Lambda \big((x_0,y_0)\big)\Big|\geq r \, ,
\]
hence \reff{app 1} and the proof of Lemma \ref{app prop} is complete. \hspace*{\fill}$\Box$\\


\bibliography{Mu_Surface}
\bibliographystyle{unsrt}

\end{document}